\documentclass[12pt]{amsart}
\usepackage{graphicx} 

\usepackage[margin=1in]{geometry}

\usepackage{amsfonts}
\usepackage{amsthm}
\usepackage{amsmath}
\usepackage{amssymb}
\usepackage{amscd} 
\usepackage{enumitem}
\usepackage{algorithm2e}
\RestyleAlgo{ruled}
\usepackage{mathtools}
\usepackage[normalem]{ulem}
\usepackage[dvipsnames]{xcolor}

\usepackage{blkarray}
\usepackage{graphicx}
\usepackage{xcolor}
\usepackage{hyperref}
\usepackage{tikz}

\usepackage{multirow}

\usepackage{verbatim}
\usepackage{url}

\newtheorem{theorem}{Theorem}[section]

\newtheorem*{thm*}{Theorem}

\theoremstyle{definition}
\newtheorem{definition}[theorem]{Definition}

\theoremstyle{remark}

\title{Suspense and surprise in matroid bingo}
\author{Matthew Baker}
\email{mbaker@math.gatech.edu}
\address{School of Mathematics, Georgia Institute of Technology, Atlanta, USA}

\date{\today}

\begin{document}

\begin{abstract}
We present the game of Matroid Bingo and prove that it gets more suspenseful the longer the game goes on.
\end{abstract}

\maketitle

Consider the following minimalist version of the game of bingo.
Each of the $m$ players gets a card on which is printed a subset of $E := \{ 1,\ldots, n \}$.
The bingo caller has a box filled with $n$ balls, on which the numbers from $1$ to $n$ inclusive are written (without duplications).
The caller repeatedly draws a uniformly random ball from among those remaining in the box, without replacement, and reads out its number.
Each player who has that number on their card crosses it out.
The caller continues this process, selecting balls at random until someone calls out ``Bingo!'', signifying that they have crossed out all the numbers on their card.

The company printing the cards wants the game to be mathematically elegant, so it imposes the following three rules:

\begin{enumerate}
    \item No bingo card is empty. In other words, no one shouts ``Bingo!'' before the game even starts.
    \item Every player has a nonzero chance of winning. In other words, no card is a proper subset of another one.
    \item There can never be a tie. Equivalently, for every possible ordering of the balls, there is a unique card that is completed first.
\end{enumerate}

\begin{figure}[h]
\centering

\begin{tikzpicture}[
  card/.style={
    draw=black!75,
    line width=0.45pt,
    rounded corners=1.5pt,
    fill=blue!32,
    minimum width=2.05cm,
    minimum height=1.18cm,
    inner sep=0pt,
    align=center,
    font=\sffamily\bfseries\large,
    text=white
  }
]
\def\dsep{\hspace{0.24em}}

\node[card] at (0,1.55) {1\dsep 2\dsep 7};
\node[card] at (2.55,1.55) {1\dsep 5\dsep 9};
\node[card] at (5.10,1.55) {2\dsep 5\dsep 7\dsep 9};

\node[card] at (0,0) {3\dsep 4\dsep 6\\[0.12em]8\dsep 9};
\node[card] at (2.55,0) {1\dsep 3\dsep 4\\[0.12em]5\dsep 6\dsep 8};
\node[card] at (5.10,0) {2\dsep 3\dsep 4\dsep 5\\[0.12em]6\dsep 7\dsep 8};
\end{tikzpicture}

\caption{A mathematical bingo deck: each of the six cards is a subset of
$E=\{1,2,3,4,5,6,7,8,9\}$.}
\label{fig:cards}
\end{figure}
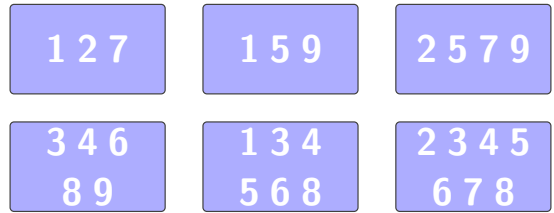

One can check that the deck of cards shown in Figure~\ref{fig:cards} satisfies all three rules.
Call such a deck of cards a \emph{mathematical bingo deck}. 
Our first surprise is that mathematical bingo decks are the same thing as \emph{matroids}.

\subsection*{Three examples}

Before defining what a matroid actually is, let us look at three examples. (Note that we previously identified $E$ with $\{1,\ldots, n\}$, but now it is convenient to allow more general labels for the elements of $E$.)

\begin{itemize}
\item In \emph{graphic bingo}, the balls are the edges of a connected finite graph $G$ and the cards are the cycles of $G$. As the caller chooses edges at random, we imagine adding those edges to the graph one at a time. The first player to complete a cycle on their card shouts ``Bingo!'' 
The game continues exactly when the set of edges that have been called so far forms a forest.\footnote{Note that for the complete graph $K_N$, revealing the edges in a random order is the Erd\H{o}s--R\'enyi random graph process, and the game clock is governed by the classical first-cycle problem, about which much is known (see, for example, \cite{FKP}).}

\item In \emph{cographic bingo}, we again start with a finite graph $G$, but now the cards are the bonds of $G$: minimal sets of edges whose removal disconnects the graph. If we assume that $G$ is connected and think of called edges as being \emph{removed} from the graph, rather than added, then the game ends as soon as connectivity is broken. Thus, graphic bingo asks when a random edge process first creates a cycle, while cographic bingo asks when random edge deletion first disconnects the graph.

\item In \emph{vector bingo}, the balls are distinct non-zero vectors $v_1,\ldots,v_n$ in a finite-dimensional vector space over some field. The cards are the minimal linearly dependent subsets, and the game ends whenever some set of called-out vectors is linearly dependent.
\end{itemize}

Matroids were invented to capture the common structure behind examples like these.

\subsection*{Matroids as mathematical bingo decks}

With those examples in mind, here is a precise definition of matroids.

\begin{definition}
A \emph{matroid} $M$ on a finite set $E$ is a collection $\mathcal C$ of subsets of $E$, called \emph{circuits}, satisfying three axioms:
\begin{itemize}
\item[(C1)] $\emptyset\notin\mathcal C.$
\item[(C2)] $\text{ If }C_1,C_2\in\mathcal C\text{ and }C_1\subseteq C_2,\text{ then }C_1=C_2.$
\item[(C3)] $\text{ If }C_1,C_2\in\mathcal C,\ C_1\ne C_2,\text{ and }e\in C_1\cap C_2,
\text{there exists }C_3\in\mathcal C\text{ such that }C_3\subseteq (C_1\cup C_2) \setminus \{e\}.$
\end{itemize}
\end{definition}

The first two axioms are precisely the first two bingo rules: the empty set is not a card, and no card contains another.
The third axiom (called the \emph{circuit elimination} axiom) is the no-tie rule in disguise. 
Indeed, suppose $C_1$ and $C_2$ are two different cards sharing a label $e$. Imagine that the caller first calls every label in
$(C_1\cup C_2) \setminus \{e\}$ and then calls $e$. If no other card had already been completed, then calling $e$ would complete both $C_1$ and $C_2$ at the same moment, creating a tie.
So the no-tie rule says that before $e$ is called, some other card $C_3$ must already have been completed. In other words,
$C_3\subseteq (C_1\cup C_2) \setminus \{e\}$,
which is exactly axiom (C3). 
Conversely, suppose that the first bingo is a tie between distinct cards $C_1$ and $C_2$, both completed when $e$ is called. 
Axiom (C3) produces a circuit $C_3\subseteq (C_1\cup C_2) \setminus \{e\}$.
But every element of $C_3$ was called before $e$, so $C_3$ was already complete --- a contradiction.

In other words, a mathematical bingo deck is the same thing as the collection of circuits of a matroid.
Matroid theorists call this a \emph{cryptomorphism}: a different-looking axiom system which nevertheless defines the same object. 
Matroids have many cryptomorphic definitions---in terms of independent sets, bases, rank functions, flats, closure operators, and more. Mathematical bingo (which from now on we will refer to as Matroid Bingo) provides a particularly intuitive cryptomorphism.

\subsection*{The game clock}

Now let us ask the question:

\medskip
\begin{center}
\boxed{\text{How long should we expect the game to last?}}
\end{center}

\medskip

Let $M$ be a matroid on $E$, with $|E|=n$, and suppose the bingo cards are the circuits of $M$.\footnote{We assume there is at least one player, since otherwise the game is rather pointless (although this does technically correspond to a matroid).}

A set $I\subseteq E$ in a matroid $M$ is called \emph{independent} if it contains no circuit. In bingo language, this means that if only the labels in $I$ have been called so far, then no one has won the game yet. The largest possible size of an independent set is called the \emph{rank} of $M$, denoted $r(M)$ (or just $r$ when the context is clear).
Our assumption that there is at least one circuit is equivalent to the statement that $r(M)<n$.

Let
\[
        i_t(M):=\#\{I\subseteq E:\ |I|=t,\ I\text{ is independent}\},
\]
and let $T$ be the round on which the first bingo occurs. 
After $t$ calls, the called set is a uniformly random $t$-element subset of $E$. The game has not ended by then exactly when that subset is independent. Therefore
\[
        \Pr(T>t)=\frac{i_t(M)}{\binom nt}.
        \tag{1}
\]

Thus the expected length of the game is
\[
\begin{aligned}
        \mathbb E[T]
        &=\sum_{t\ge0} t \cdot \Pr(T=t) \\
        &=\sum_{t\ge0}\Pr(T>t) \\
        &=\sum_{t=0}^{r(M)}\frac{i_t(M)}{\binom nt}. \\
\end{aligned}        
\]
So the whole game-clock distribution depends only on $n$, $r$, and the independent-set numbers
\[
        i_0(M),i_1(M),\ldots,i_{r}(M).
\]
These numbers are among the most widely studied numerical invariants of a matroid.


For $t\geq 1$, define the \emph{suspense level} just before the $t^{\rm th}$ call to be
\[
        h_t:=\Pr(T=t\mid T>t-1).
\]
Thus $h_t$ is the probability that the $t^{\rm th}$ call ends the game, conditional on the game having survived the first $t-1$ calls.
(Probabilists call this a \emph{hazard rate}, but ``suspense'' feels better suited to a non-lethal game like bingo.)
Using (1), we find that
\[
        h_t
        =1-\frac{\Pr(T>t)}{\Pr(T>t-1)}
        =1-
        \frac{i_t(M)/\binom nt}
             {i_{t-1}(M)/\binom n{t-1}}
        \tag{2}
\]
for $1 \leq t \leq r+1$.

\begin{table}[h]
\centering
\small
\begin{tabular}{r|ccccccc}
$t$      & 1 & 2 & 3 & 4 & 5 & 6 & 7 \\
\hline
$h_t$    
         & $.000$ & $.000$ & $.024$ & $.081$ & $.186$ & $.397$ & $1$ \\[3pt]
\end{tabular}
\caption{The suspense levels for the bingo deck of
Figure~\ref{fig:cards}, which has rank $r=6$ on $n=9$ elements.
Note that the sequence $h_t$ is weakly increasing.}
\label{tab:suspense}
\end{table}

We now ask the question:

\medskip

\begin{center}
\boxed{\text{Does bingo get more suspenseful as time goes on?}}
\end{center}

\medskip

And now for another surprise: to answer this question, we must turn to a deep recent advance in matroid theory which is connected to both June Huh's 2022 Fields Medal and Shayan Oveis Gharan's 2026 Abacus Medal.

\subsection*{Combinatorial Hodge theory}

A sequence $a_0,a_1,a_2,\ldots, a_d$ of nonnegative real numbers
is called \emph{log-concave} if
\[
        a_t^2\ge a_{t-1}a_{t+1}
\]
for $t=1,\ldots,d-1$. When the $a_t$'s are all positive, this is equivalent to saying that the ratios $a_t/a_{t-1}$ are weakly decreasing.

The strong version of Mason's conjecture, which was proved independently in \cite[Theorem 4.14]{BH} and \cite[Theorem 1.2]{ALOV}, says:

\begin{theorem}[Mason's conjecture, strong form]
If $M$ is a matroid of rank $r$ on $n$ elements, then the sequence
$i_t(M)/\binom nt$ is log-concave for $t=0,\ldots,n$.
\end{theorem}

Equivalently, for all $1 \leq t \leq n-1$,
\[
\left(\frac{i_t(M)}{\binom nt}\right)^2
\ge
\left(\frac{i_{t-1}(M)}{\binom n{t-1}}\right)
\left(\frac{i_{t+1}(M)}{\binom n{t+1}}\right).
\]

The proof of the strong Mason theorem is part of one of the most exciting recent stories in algebraic combinatorics, namely the emergence of Hodge-theoretic methods in matroid theory. While the margins of this manuscript are too small to fit the proof, we can at least explain precisely what Mason's theorem says about our game of bingo.

\subsection*{The main result}

Set
\[
        S(t):=\Pr(T>t)=\frac{i_t(M)}{\binom nt}.
\]
If $r=r(M)$, then $S(t)>0$ for $0\le t\le r$ and $S(t)=0$ for all $t>r$. 
Formula (2) says
\[
        h_t=1-\frac{S(t)}{S(t-1)}
\]
for $1\le t\le r$, while $h_{r+1}=1$ (since every $(r+1)$-element set is dependent).

The strong Mason theorem tells us that $S(0),S(1),\ldots,S(r+1)$ is log-concave, and thus the ratios
$S(t)/S(t-1)$ are weakly decreasing for $1\le t\le r+1$. Hence the sequence $h_1,h_2,\ldots,h_{r+1}$ 
is weakly increasing. We have thus proved:

\begin{theorem}[Increasing Suspense Theorem]
In any legal game of mathematical bingo, conditional on the game not having ended yet, the probability that the next call ends the game never decreases. 
\end{theorem}

In other words: the longer the game lasts, the more suspense the next ball holds.

\subsection*{Why Hodge theory?}

Hodge theory originated in algebraic geometry and topology, where it gives powerful constraints on the cohomology of smooth projective varieties. Among its central consequences are the Hard Lefschetz theorem and the Hodge--Riemann relations, which force certain numerical sequences to be log-concave.
Over the last decade, June Huh and his collaborators discovered that matroids possess analogous Hodge-theoretic structures. For example, Adiprasito, Huh, and Katz constructed a Chow ring for matroids and proved Hodge--Riemann relations for it, settling several long-standing log-concavity conjectures \cite{AHK}. Later developments, including Br\"and\'en and Huh's theory of Lorentzian polynomials \cite{BH} and the closely related work of Anari, Liu, Oveis Gharan, and Vinzant on completely log-concave polynomials \cite{ALOV}, led to the strong form of Mason's conjecture used above.


\subsection*{Individual cards}

The Increasing Suspense Theorem is purely about the game clock and says nothing about which player is most likely to win.
That question is a bit subtler. For example, if one card has fewer numbers on it than another, should you always choose the smaller card? 

Our final surprise is that the answer to this question is typically yes, but occasionally no.
In joint work with Hope Dobbelaere, Brennan Fullmer, and Patrik Gajdo\v{s}, we studied the winning probabilities of individual cards in Matroid Bingo \cite{BDFG}. We found that longer cards can sometimes be more likely to win than shorter ones; we called this a \emph{monotonicity violation}. Such violations exist, but they are surprisingly rare:
among the 385,360 isomorphism classes of matroids on at most nine elements that have at least one circuit --- and hence define a nonempty bingo game --- exactly 11 exhibit a monotonicity violation.
One of these monotonicity violations is depicted in Figure~\ref{fig:cards}: if given a choice between $\{ 2, 5, 7, 9 \}$ and $\{ 3, 4, 6, 8, 9 \}$, it turns out that you are better off choosing the latter.
Indeed,
\[
{\rm Pr}\left(34689 \text{ wins}\right) = \frac{37}{252} > \frac{13}{90} = {\rm Pr}\left(2579 \text{ wins}\right)
\]
since $\frac{37}{252} \approx 0.1468$ and $\frac{13}{90} \approx 0.1444$.

Our paper also proves that, for any fixed card $C$, the probability that $C$ wins in round $t$ is log-concave and has no internal zeroes, and in particular it is unimodal. This also requires an application of the strong Mason theorem.

\subsection*{Credits}

Matroid Bingo was introduced by Dima Fon-Der-Flaass and is mentioned, for example, in a book by Borovik, Gelfand, and White \cite{BGW}. It was further explored in \cite{BDFG}.
The Increasing Suspense Theorem is new to this note.

\subsection*{Acknowledgments}

We thank Hope Dobbelaere, Brennan Fullmer, and Patrik Gajdo\v{s} for helpful discussions.
The author was supported by NSF grant DMS-2154224 and a Simons Fellowship in Mathematics.

\subsection*{Statement of AI usage}
We thank GPT-5.5 Pro for assistance with polishing the exposition in this note, and GPT-5.6 Pro for creating the TikZ figure and proofreading the manuscript.

\end{document}